\documentclass[12pt,reqno]{amsart}
\usepackage{enumerate, latexsym, amsmath, amsfonts, amssymb, amsthm, color}
\usepackage{psfrag,graphicx, colordvi}
\newtheorem{theo}{Theorem}

\newtheorem{coro}[theo]{Corollary}

\newtheorem{lem}[theo]{Lemma}
\def\pmod #1{\ ({\rm{mod}}\ #1)}
\def\Z{\Bbb Z}
\def\N{\Bbb N}

\def\l{\left}
\def\r{\right}
\def\bg{\bigg}
\def\({\bg(}
\def\){\bg)}
\def\t{\text}
\def\f{\frac}

\def\ls{\leqslant}

\def\bi{\binom}

\def\eq{\equiv}

\def\Proof{\noindent{\it Proof}}

\theoremstyle{plain}

\theoremstyle{definition}

\theoremstyle{remark}

\makeatletter \@addtoreset{equation}{section}
\@addtoreset{theo}{section} \makeatother

\def\qed{\hfill \rule{4pt}{7pt}}

\begin{document}
\baselineskip=17pt
\hbox{Int. J. Number Theory 14(2018), no.\,1, 143--165.}
\medskip

\title
[Telescoping Method and
Congruences for
Double Sums]
{Telescoping Method and \\Congruences for
Double Sums}

\author
[Yan-Ping Mu and Zhi-Wei Sun] {Yan-Ping Mu and Zhi-Wei Sun}

 \address {(Yan-Ping Mu) College of Science, Tianjin University of Technology,
 Tianjin 300384, P. R. China}

\email{{\tt yanping.mu@gmail.com}}

\address {(Zhi-Wei Sun) Department of Mathematics, Nanjing
University, Nanjing 210093, People's Republic of China}
\email{{\tt zwsun@nju.edu.cn}}

\keywords{Telescoping method, double sum, combinatorial sequence, congruence.
\newline \indent 2010 {\it Mathematics Subject Classification}. Primary 11B75, 33F10; Secondary 05A10, 11A07.}
\thanks{The second author is the corresponding author. Both authors are supported by the National Natural Science Foundation of China (grants 11471244 and 11571162, respectively)}

\begin{abstract}
In recent years, Z.-W. Sun proposed several sophisticated conjectures on congruences for finite sums with terms involving combinatorial sequences such as
central trinomial coefficients, Domb numbers and Franel numbers. These sums are double summations of hypergeometric terms.
 Using the telescoping method and certain mathematical software packages, we transform such a double summation into a single sum. With this new approach, we confirm several open conjectures of Sun.
\end{abstract}

\maketitle

\section{Introduction}

For a sequence of integers $a_0,a_1,a_2,\ldots$, it is interesting to see whether the arithmetic mean
$$\f{a_0+a_1+\cdots+a_{n-1}}{n}$$
is integral for any $n\in\Z^+=\{1,2,3,\ldots\}$. It is easy to see that
$$\f1n\sum_{k=0}^{n-1}(2k+1)(-1)^k=(-1)^{n-1}\in\Z\ \ \mbox{and}\ \ \f1{n^2}\sum_{k=0}^{n-1}(2k+1)=1\in\Z.$$
For the sequence of Ap\'ery numbers
$$A_k:=\sum_{l=0}^k\bi kl^2\bi{k+l}l^2=\sum_{l=0}^k\bi{k+l}{2l}^2\bi{2l}l^2\ \ (k\in\N=\{0,1,2,\ldots\}),$$
Z.-W. Sun \cite{S12} proved the congruence $\sum_{k=0}^{n-1}(2k+1)A_k\eq0\pmod n$, and conjectured that
$\sum_{k=0}^{n-1}(2k+1)(-1)^kA_k\eq0\pmod n$, which was confirmed by V.J.W. Guo and J. Zeng \cite{GZa}.
For the sequence of Franel numbers
$$f_k:=\sum_{l=0}^k\bi kl^3\ (k=0,1,2,\ldots)$$
(for which J. Franel \cite{Fra95} found a recurrence),
Z.-W. Sun \cite{S13b} conjectured that
$$\sum_{k=0}^{n-1}(3k+2)(-1)^kf_k\eq0\pmod{2n^2},$$
which was confirmed by Guo \cite{Guo}. However, some of such conjectural congruences still remain open.
For example, Sun \cite[Conjecture 5.1]{S13c} conjectured that
$$\sum_{k=0}^{n-1}(5k+4)D_k\eq0\pmod{4n},$$
where $D_0,D_1,D_2,\ldots$ are the Domb numbers defined by
$$D_k:=\sum_{l=0}^k\bi kl^2\bi{2l}l\bi{2(k-l)}{k-l}\ \ (k\in\N).$$
For various combinatorial interpretations of the Domb numbers, see L.B. Richmond and J. Shallit \cite{RS}, and N.J.A Sloane \cite{S}.

In this paper we study some conjectures of Sun \cite{S14, Sun16} on congruences for sums of the form
\begin{equation}\label{genf}
S_n = \sum_{k=0}^{n-1} \sum_{l=0}^k F(k,l)\ \ (n\in\Z^+),
\end{equation}
where $F(k,l)$ is a bivariate hypergeometric term of $k$ and $l$.
To confirm such conjectures, we try to search for two hypergeometric terms $G_1(k,l)$ and $G_2(k,l)$ such that
\begin{equation}\label{bi-tel}
F(k,l) = \Delta_k (G_1(k,l)) + \Delta_l (G_2(k,l)),
\end{equation}
where
\[G_1(k,l)=R_1(k,l)F(k,l)\quad \mbox{and}\quad   G_2(k,l)=R_2(k,l)F(k,l)\]
with $R_1(k,l)$ and $R_2(k,l)$ being rational functions,
and
$$\Delta_k (G_1(k,l)) = G_1(k+1,l)-G_1(k,l)$$
and $$\Delta_l (G_2(k,l)) = G_2(k,l+1)-G_2(k,l).$$
Though the denominator of the rational function $R_1(k,l)$ or $R_2(k,l)$ might be zero for some $0\ls l\ls k\ls n-1$,
the functions $G_1(k,l)$ and $G_2(k,l)$ we obtain are essentially well defined for $0\ls l\ls k\ls n-1$ and hence (\ref{bi-tel}) can be verified directly.
Once we have $G_1(k,l)$ and $G_2(k,l)$ in hand, the sum $S_n$ can be transformed to a single sum
\begin{equation}\label{Sn}
S_n = \sum_{l=0}^{n-1} \big( G_1(n,l)-G_1(l,l) \big)  + \sum_{k=0}^{n-1} \big( G_2(k,k+1)-G_2(k,0) \big) .
\end{equation}

We can use the {\tt Maple} package {\it DoubleSum} given by Chen-Hou-Mu \cite{Chen06} or the {\tt Mathematica} package {\it HolonomicFunctions} given by C. Koutschan \cite{Kou10} to compute $G_1(k,l)$ and $G_2(k,l)$. Usually, the two hypergeometric terms $G_1(k,l)$ and $G_2(k,l)$ obtained directly by these two packages are not optimal. We will remove some factors from the denominators of the resulting $R_1$ and $R_2$, and use the package {\it DoubleSum} or the package {\it MultiSum} \cite{Lyo02} to compute a suitable pair $(G_1(k,l), G_2(k,l))$.  Once we get the simplification \eqref{Sn} for $S_n$, it would be be convenient to deduce Sun's conjectural congruences for $S_n$.
Using this powerful method, we confirm several sophisticated open conjectures of Sun.

The second author \cite{S13a, S13b} established several supercongruences on sums of Franel numbers.
Our first theorem confirms a conjecture of him involving Franel numbers.

\begin{theo}[\mbox{\cite[Conjecture 4.1(i)]{Sun16}}] \label{s1}
For any integer $n>1$, we have
\begin{equation}\label{con1}
\sum_{k=0}^{n-1} (9k^2 + 5k) (-1)^k f_k \equiv 0 \pmod{n^2(n-1) }.
\end{equation}
Moreover, for any odd prime $p$ we have
\begin{equation}\label{f-cong}
\sum_{k=0}^{p-1} (9k^2 + 5k) (-1)^k f_k \equiv  3p^2(p-1)-16p^3q_p(2) \pmod{p^4},
\end{equation}
where $q_p(2)$ denotes the Fermat quotient $(2^{p-1}-1)/p$.
\end{theo}

In \cite{Sun16} Sun introduced the polynomials
$$g_n(x):=\sum_{k=0}^n\bi nk^2\bi{2k}kx^k\ (n=0,1,2,\ldots)$$
and obtained some related congruences. Such polynomials are interesting since it is closely related to the Ap\'ery polynomials
$$A_n(x):=\sum_{k=0}^n\bi nk^2\bi{n+k}k^2x^k=\sum_{k=0}^n\bi{n+k}{2k}^2\bi{2k}k^2x^k\ (n\in\N);$$
in fact, Sun \cite[(2.8)]{Sun16} showed that
$$A_n(x)=\sum_{k=0}^n\bi nk\bi{n+k}k(-1)^{n-k}g_k(x).$$

\begin{theo}[\mbox{\cite[Conjecture 4.1(ii)]{Sun16}}] \label{s2}
For any $n\in\Z^+$, the number
\[
\frac{1}{n^2} \sum_{k=0}^{n-1} (8k^2+12k+5) g_k(-1)
\]
is an odd integer. Moreover, for any prime $p$ we have
\begin{equation}\label{mp3}
\sum_{k=0}^{p-1} (8k^2+12k+5) g_k(-1) \equiv 3p^2 \pmod{p^3}.
\end{equation}
\end{theo}
{\noindent \it Remark}.  Guo, Mao and Pan \cite{Guo15} called $g_n(x)\ (n\in\N)$ {\it Sun polynomials}, and they used the Zeilberger algorithm (cf. \cite[pp.\,101-119]{PWZ}) to prove that
\[
 \sum_{k=0}^{n-1} (8k^2+12k+5) g_k(-1) \equiv 0 \pmod{n}\ \ \t{for all}\ n\in\Z^+,
\]
 which is the first progress on \cite[Conjecture 4.1(ii)]{Sun16}.

\begin{theo}[\mbox{\cite[Conjecture 4.3]{Sun16}}] \label{s3}
For $k=0,1,2,\ldots$ define
$$F_k := \sum_{l=0}^k {k \choose l}^3 (-8)^l.$$
Then, for any positive integer $n$, the number
\[\frac{1}{n} \sum_{k=0}^{n-1} (6k+5) (-1)^k F_k\]
is an odd integer.
\end{theo}

For $n\in\N$, the {\it central trinomial coefficient} $T_n$ denotes the coefficient of $x^n$ in the expansion of $(x^2+x+1)^n$.
$T_n$ has many combinatorial interpretations. For example, $T_n$ counts the number of permutations of $n$ symbols, each $-1$, $0$, or $1$, which sum to $0$.
For $b,c\in\Z$, Sun \cite{S14} defined $T_n(b,c)$ as the coefficient of $x^n$ in the expansion of $(x^2+bx+c)^n$, and called those $T_n(b,c)\ (n\in\N)$ {\it generalized central trinomial coefficients}.

\begin{theo}[\mbox{\cite[Conjecture 5.2]{S14}}]
Let $b,c \in \mathbb{Z}$. For any $n\in\Z^+$, we have
\begin{equation}
\sum_{k=0}^{n-1} (8ck+4c+b)T_k(b,c^2)^2 (b-2c)^{2(n-1-k)} \equiv 0 \pmod{n}.
\end{equation}
If $p$ is an odd prime not dividing $b(b-2c)$, then
\begin{equation}
\sum_{k=0}^{p-1} (8ck+4c+b) \frac{T_k(b,c^2)^2}{(b-2c)^{2k}}\equiv p(b+2c) \left( \frac{b^2-4c^2}{p} \right) \pmod{p^2},
\end{equation}
where $(\f{\cdot}p)$ denotes the Legendre symbol.
\end{theo}

In the case $b=c=1$, Theorem 1.4 yields the following consequence.

\begin{coro}[\mbox{\cite[Conjecture 1.1(i)]{S14}}]
For any positive integer $n$, we have
\begin{equation}
\sum_{k=0}^{n-1} (8k+5)T_k^2 \equiv 0 \pmod{n}.
\end{equation}
If $p$ is a prime, then
\begin{equation}
\sum_{k=0}^{p-1} (8k+5)T_k^2 \equiv 3p \left( \frac{p}{3} \right) \pmod{p^2}.
\end{equation}
\end{coro}

Our next theorem confirms a conjecture of Sun involving the Domb numbers.

\begin{theo}[\mbox{\cite[Conjecture 5.15(ii)]{S11a}}] \label{s4}
For any integer $n > 1$, the number
\[
a_n = \frac{1}{2n^3(n-1)} \sum_{k=0}^{n-1} (3k^2+k)D_k 16^{n-1-k}
\]
is an integer. Moreover, $a_n$ is odd if and only if $n$ is a power of two.
Also, for any prime $p>3$ we have the supercongruence
\begin{equation}\label{cd}\sum_{k=0}^{p-1}\f{3k^2+k}{16^k}D_k\eq-4p^4q_p(2)\pmod{p^5}.\end{equation}
\end{theo}

Our last theorem is about some conjectures of the second author
posed in a preprint form {\tt arxiv:1407.0967v6} of \cite{Sun16}.

\begin{theo} \label{s7}
{\rm (i)} For any integer $n>1$, we have
\begin{equation}\label{f3}\sum_{k=0}^{n-1}(12k^4+25k^3+21k^2+6k)(-1)^kf_k\eq0\pmod{4n^2(n-1)}.
\end{equation}
If $p$ is an odd prime, then
\begin{equation}\label{f4}\sum_{k=0}^{p-1}(12k^4+25k^3+21k^2+6k)(-1)^kf_k\eq-4p^3\pmod{p^4}.
\end{equation}

{\rm (ii)} For any $n\in\Z^+$ we have
\begin{equation}\label{g3}\sum_{k=0}^{n-1}(12k^3+34k^2+30k+9)g_k\eq0\pmod{3n^3},
\end{equation}
where $g_k:=g_k(1)=\sum_{l=0}^k\bi kl^2\bi{2l}l$.
If $p$ is an odd prime, then
\begin{equation}\label{g4}\sum_{k=0}^{p-1}(12k^3+34k^2+30k+9)g_k\eq\f{3p^2}2\l(1+3\l(\f p3\r)\r)\pmod{p^4}.
\end{equation}

{\rm (iii)} For any $n\in\Z^+$ we have
\begin{equation}\label{a4}\sum_{k=0}^{n-1}(18k^5+45k^4+46k^3+24k^2+7k+1)(-1)^kA_k\eq0\pmod{n^4}.\end{equation}
Also, for any prime $p>3$ we have
\begin{equation}\label{a7}\sum_{k=0}^{p-1}(18k^5+45k^4+46k^3+24k^2+7k+1)(-1)^kA_k\eq-2p^4+3p^5\pmod{p^7}.
\end{equation}
\end{theo}
\noindent{\it Remark}. The second author even conjectured that
\begin{equation}\label{a10}
\begin{aligned}&\sum_{k=0}^{p-1}(18k^5+45k^4+46k^3+24k^2+7k+1)(-1)^kA_k
\\\eq& -2p^4+3p^5+(6p-8)p^5H_{p-1}-\f{12}5p^9B_{p-5}\pmod{p^{10}}
\end{aligned}
\end{equation}
for any prime $p>5$, where $H_n$ denotes the $n$-th harmonic number $\sum_{0<k\ls n}1/k$, and $B_0,B_1,B_2,\ldots$ are the Bernoulli numbers.

\section{Proof of Theorem~\ref{s1}}

In this section, we illustrate in detail how to prove Theorem~\ref{s1} by the telescoping method.

\medskip
\noindent {\it Proof of Theorem} \ref{s1}. By an identity of V. Strehl \cite{St}, for any $k\in\N$ we have
\begin{equation}\label{f-iden}
f_k = \sum_{l=0}^k {k \choose l}^2 {2l \choose k}.
\end{equation}
Therefore, for any integer $n>1$, the left hand side of the congruence \eqref{con1} can be written as
\[
\sum_{k=0}^{n-1} \sum_{l=0}^k F(k,l),
\]
where
\[
F(k,l) = (9k^2+5k) (-1)^k {k \choose l}^2 {2l \choose k}.
\]
By the command {\tt GetDen(F, k, l)} of the package {\it DoubleSum}, we obtain the estimated denominators of $R_1(k,l)$ and $R_2(k,l)$:
\begin{align*}
d_1(k,l) & = k(9k+5)(l+1)(-2l+k-1) (-2l-2+k),\\
d_2(k,l) & = k(9k+5)(k-l+1)^2.
\end{align*}
We remove the factor $k(-2l-2+k)$ from $d_1(k,l)$, and the factor $(k-l+1)^2$ from $d_2(k,l)$, and then use the command
\begin{verbatim}
  SolveR(F, n, k, l, (9k+5)(l+1)(-2l+k-1), k(9k+5), 0, "s")
\end{verbatim}
to obtain the rational functions
$$R_1(k,l)  = \frac{(k-l)^2 (3k+4l-3)}{(9k+5)(l+1)(-2l+k-1)}$$
and
$$ R_2(k,l)  = \frac{(-2l+k)(3k-l-2)}{k(9k+5)}.$$
Therefore
\[
F(k,l) = \Delta_k (G_1(k,l)) + \Delta_l (G_2(k,l)),
\]
where
\begin{align*}
G_1(k,l) & = (-1)^{k-1}k^2(3k+4l-3) \frac{1}{l+1}{2l \choose l} {k-1 \choose l}{l \choose k-1-l}, \\
G_2(k,l) & = (-1)^{k-1}(2l-k)(3k-l-2) {k \choose l}^2 {2l \choose k}.
\end{align*}
Noting that
\[
G_1(l,l)=G_2(k,k+1)=G_2(k,0) = 0,
\]
we get
\begin{align*}
S_n  =& \sum_{l=0}^{n-1} G_1(n,l) \\
 =& (-1)^{n-1} n^2\sum_{l=0}^{n-1} (3n+4l-3) \frac{1}{l+1}{2l \choose l} {n-1 \choose l}{l \choose n-1-l} \\
 =& (-1)^{n-1}3 n^2(n-1)   \sum_{l=0}^{n-1} \frac{1}{l+1}{2l \choose l} {n-1 \choose l}{l \choose n-1-l} \\
&+ (-1)^{n-1}4 n^2(n-1) \sum_{l=1}^{n-1}  \frac{1}{l+1}{2l \choose l} {n-2 \choose l-1}{l \choose n-1-l} .
\end{align*}
Since $C_l={2l \choose l}/(l+1)$ is the $l$-th Catalan number which is an integer, we derive immediately that $S_n$ is divisible by $n^2(n-1)$.

Now let $p$ be an odd prime. By the above,
\begin{equation}\label{con1-p}
S_p=p^2  \sum_{l=(p-1)/2}^{p-1} \f{3p-3+4l}{l+1} {p-1 \choose l}\bi{2l}{l}\bi{l}{p-1-l}.
\end{equation}
For $l=(p+1)/2,\ldots,p-2$, we clearly have $p\mid\bi{2(l+1)}{l+1}$ and hence
\begin{align*}\bi{2l}l\bi l{p-1-l}=&\f{l+1}{2l+1}\bi{2l+1}{l}(-1)^l\bi{-l+(p-1-l)-1}{p-1-l}
\\=&\f{(-1)^l(l+1)}{2(2l+1)}\bi{2(l+1)}{l+1}\bi{2(p-1-l)-p}{p-1-l}
\\\eq&\f{(-1)^l(l+1)}{2(2l+1)}\bi{2(l+1)}{l+1}\bi{2(p-(l+1))}{p-(l+1)}
\\\eq&\f{(-1)^lp}{2l+1}\pmod{p^2}
\end{align*}
with the help of \cite[Lemma 2.1]{S11a}. Therefore
\begin{align*}\f {S_p}{p^2}\eq&\sum_{(p-1)/2<l<p-1}\f{3p-3+4l}{l+1} {p-1 \choose l}\f{(-1)^lp}{2l+1}
\\&+\f{3p-3+4(p-1)/2}{(p+1)/2}\bi{p-1}{(p-1)/2}^2+\f{3p-3+4(p-1)}p\bi{2(p-1)}{p-1}
\\\eq&p\sum_{(p-1)/2<l<p-1}\f{4l-3}{(l+1)(2l+1)}+\f{10(p-1)}{p+1}(1-pH_{(p-1)/2})^2
\\&+\f{7(p-1)}{2p-1}\bi{2p-1}{p-1}\pmod{p^2}
\end{align*}
by noting that
$$\bi{p-1}{(p-1)/2}(-1)^{(p-1)/2}=\prod_{j=1}^{(p-1)/2}\l(1-\f pj\r)\eq 1-pH_{(p-1)/2}\pmod{p^2}.$$
Clearly,
$$H_{p-1}=\sum_{k=1}^{(p-1)/2}\l(\f1k+\f1{p-k}\r)\eq0\pmod{p}$$
and
$$\bi{2p-1}{p-1}=\prod_{k=1}^{p-1}\l(1+\f pk\r)\eq 1+pH_{p-1}\eq1\pmod{p^2}.$$
Thus, we deduce that
\begin{align*}\f {S_p}{p^2}\eq&p\sum_{(p-1)/2<l<p-1}\l(\f7{l+1}-\f{10}{2l+1}\r)
\\&-10(p-1)^2\l(1-2pH_{(p-1)/2}\r)-7(p-1)(2p+1)
\\\eq&7p\l(H_{p-1}-H_{(p-1)/2}-\f2{p+1}\r)-10p\sum_{j=1}^{(p-3)/2}\f1{p+2j}
\\&+10(2p-1)(1-2pH_{(p-1)/2})+7(p+1)
\\\eq&-7pH_{(p-1)/2}-14p-10p\l(\f{H_{(p-1)/2}}2-\f1{p-1}\r)
\\&+10(2p-1+2pH_{(p-1)/2})+7(p+1)
\\\eq&8pH_{(p-1)/2}-24p+27p-3
\pmod{p^2}.
\end{align*}
As
$$-\sum_{k=1}^{(p-1)/2}\f p{2k}\eq\sum_{k=1}^{(p-1)/2}\f{p}{2k}\bi{p-1}{2k-1}=\sum_{k=1}^{(p-1)/2}\bi p{2k}=2^{p-1}-1\pmod{p^2},$$
we have $pH_{(p-1)/2}\eq-2(2^{p-1}-1)\pmod {p^2}$, i.e.,
\begin{equation}\label{hh}H_{(p-1)/2}\eq-2q_p(2)\pmod p.\end{equation}
 So, we finally obtain the congruence
$$\frac{S_p}{p^2}\eq8p(-2q_p(2))+3p-3\pmod{p^2}$$
which gives (\ref{f-cong}).
This concludes our proof of Theorem~\ref{s1}. \qed

\section{Proofs of Theorems~\ref{s2} and \ref{s3}}

In this section, we prove Theorems ~\ref{s2} and \ref{s3}. The method is similar to the one used in Section 2. We will give related $R_1(k,l)$ and $R_2(k,l)$ directly.

\medskip
\noindent{\it Proof of Theorem} \ref{s2}.
Let
\[
F(k,l) = (8k^2 + 12k + 5) {k \choose l}^2 {2l \choose l} (-1)^l
\]
be the summand and
\[
S_n = \sum_{k=0}^{n-1} \sum_{l=0}^k F(k,l)\qquad\ \mbox{with}\ n\in\Z^+.
\]
We find that
\begin{align*}
R_1(k,l)  =  \frac{ (k-l)^2(2k+5l+3)}{(8k^2+12k+5)(l+1)} \ \mbox{and}\
R_2(k,l)  = - \frac{l^2}{8k^2+12k+5}.
\end{align*}
Thus we derive that
\[
S_n = \sum_{l=0}^{n-1} R_1(n,l)F(n,l) = n^2 \sum_{l=0}^{n-1} (2n+5l+3)(-1)^l  {n-1 \choose l}^2 \frac{1}{l+1} {2l \choose  l}.
\]
Noting that $C_l={2l \choose  l}/(l+1)\in\Z$, we have
\[
S_n\equiv n^2\sum_{l=0}^{n-1} (-1)^l  {n-1 \choose l}^2 {2l \choose  l}
\eq n^2 \pmod{2n^2}
\]
since$\bi{2l}l=2\bi{2l-1}{l-1}$ for $l=1,2,3,\ldots$.
Therefore, $S_n/n^2$ is an odd integer.

It is straightforward to check that
$S_2 = -20 \equiv 3 \times 2^2 \pmod{2^3}.$
So ~\eqref{mp3} holds for $p=2$.
Now let $p$ be an odd prime. Since
\[
{p-1 \choose l} \equiv (-1)^l \pmod{p}\quad \ \mbox{for all}\ l=0,1,\ldots,p-1,
\]
we have
\begin{align*}
\frac{S_p}{p^2} & \equiv  \sum_{l=0}^{p-1} (5l+3) \frac{(-1)^l}{l+1} {2l \choose  l}  =  5 \sum_{l=0}^{p-1} (-1)^l {2l \choose  l}  - 2 \sum_{l=0}^{p-1} (-1)^l C_l \pmod{p}.
\end{align*}
Clearly,
\begin{align*}
\sum_{l=0}^{p-1} (-1)^l {2l \choose  l} = &\sum_{l=0}^{p-1}\bi{-1/2}l4^l
\\\eq&\sum_{l=0}^{p-1}\bi{(p-1)/2}l4^l=5^{(p-1)/2}\eq\l(\f 5p\r)=\l(\f p5\r)\pmod p
\end{align*}
and
$$\sum_{l=1}^{p-1}(-1)^lC_l\eq-\f52\l(1-\l(\f 5p\r)\r)=\f 52\l(\l(\f p5\r)-1\r)\pmod{p}$$
by \cite[Lemma 2.1]{S11b}.
We thus derive that
\[
S_p \equiv  3p^2 \pmod{p^3},
\]
completing the proof of Theorem~\ref{s2}.
\qed

\medskip
\noindent{\it Proof of Theorem} 1.3. For $n\in\Z^+$ let
\[S_n = \sum_{k=0}^{n-1}(6k+5)(-1)^kF_k = \sum_{k=0}^{n-1} \sum_{l=0}^k F(k,l)\]
with $F(k,l)=(6k+5) (-1)^k {k \choose l}^3 (-8)^l$.
We find that
$$R_1(k,l)  = - \frac{4}{3}\cdot \frac{(12 k^2-30 k l+21 l^2-32 k+46 l+25)(k-l)^3 }{(6k+5) (l+1)^3}$$
and $$R_2(k,l) = \frac{1}{3}\cdot \frac{12 k^2-12 k l+12 l^2+10 k+4 l+5}{6 k+5}.$$
We thus derive that
\begin{align*}
S_n &= \sum_{l=0}^{n-1} R_1(n,l)F(n,l) - \sum_{k=0}^{n-1} R_2(k,0)F(k,0) \\
& = - \frac{4}{3} \sum_{l=0}^{n-1} (12 n^2-30 n l+21 l^2-32 n+46 l+25) {n \choose l+1}^3 (-8)^l (-1)^n \\
& \ \ -  \frac{1}{3} \sum_{k=0}^{n-1} (12 k^2 + 10k + 5) (-1)^k.
\end{align*}
The second summation can be easily evaluated directly or by {\tt Maple}:
\[
-\frac{1}{3} \sum_{k=0}^{n-1} (12 k^2 + 10k + 5) (-1)^k = -\frac{n}{3} (-1)^n + 2 (-1)^n n^2.
\]
Note that $(l+1)\bi n{l+1}=n\bi{n-1}l\eq0\pmod n$. So we have
\begin{align*}
(-1)^nS_n\eq&-\f 43\sum_{l=0}^{n-1}(21l+25)(l+1)\bi n{l+1}^3(-8)^l
\\&-\f43n\sum_{l=0}^{n-1}\bi n{l+1}^3(-8)^l-\f n3
\\\eq&-\f 43\sum_{l=0}^{n-1}(l+1)\bi n{l+1}^3-\f 43n\sum_{j=0}^n\bi nj^3+n
\\=&-\f23\sum_{j=0}^n(j+(n-j))\bi nj^3-\f 43n\sum_{j=0}^n\bi nj^3+n
\\\eq&-2nf_n+n \eq n\pmod {2n}.
\end{align*}
This proves the desired result. \qed

\section{Proof of Theorem 1.4}

\medskip
\noindent{\it Proof of Theorem} 1.4. By \cite[Lemma 4.1]{S14}, for any $k\in\N$ we have
\[
T_k(b,c^2)^2 = \sum_{l=0}^k {k+l \choose 2l} {2l \choose l}^2 c^{2l} (b^2-4c^2)^{k-l}.
\]
Let
\[
F(k,l) = (8ck+4c+b) {k+l \choose 2l} {2l \choose l}^2 c^{2l} (b^2-4c^2)^{k-l}(b-2c)^{2(n-1-k)}
\]
be the summand and
\[
S_n = \sum_{k=0}^{n-1} \sum_{l=0}^k F(k,l)=\sum_{k=0}^{n-1}(8ck+4c+b)T_k(b,c^2)^2(b-2c)^{2(n-1-k)}
\]
for any $n\in\Z^+$.

If $b=-2c$, then $T_k(b,c^2)^2=\bi{2k}k^2c^{2k}$ for all $k\in\N$, and hence
\begin{align*}&\sum_{k=0}^{n-1}(8ck+4c+b)T_k(b,c^2)^2(b-2c)^{2(n-1-k)}
\\=&\sum_{k=0}^{n-1}2c(4k+1)\bi{2k}k^2c^{2k}(-4c)^{2(n-1-k)}
\\=&2c16^{n-1}\sum_{k=0}^{n-1}(4k+1)\f{\bi{2k}k^2}{16^k}=2cn^2\bi{2n-1}n^2
\\\eq&0\pmod{n^2}.
\end{align*}
(The identity $\sum_{k=0}^{n-1}(4k+1)\bi{2k}k^2/16^k=\bi{2n-1}n^2n^2/16^{n-1}$ can be easily proved by induction.)

Below we assume that $b+2c\not=0$. We find the rational functions
$$R_1(k,l) = -\frac{(k-l)(-b^2+8c^2+lb^2+8lc^2-2kbc+2b^2l^2-4kbcl)}{(b+2c)(8ck+4c+b)(l+1)}$$
and $$R_2(k,l) = -\frac{4bl^2}{8ck+4c+b}.$$
It follows that
\begin{align*}
S_n =& \sum_{l=0}^{n-2} R_1(n,l)F(n,l) + R_1(n,n-1)F(n,n-1) \\
=& -n \sum_{l=0}^{n-2} (-b^2+8c^2+lb^2+8lc^2-2nbc+2b^2l^2-4nbcl)\\
\\& \ \ \times c^{2l}(b+2c)^{n-l-1}(b-2c)^{n-l-2}
{n-1 \choose l}{n+l \choose l}C_l\\
& -n^2(2nb-3b-4c)c^{2n-2} {2n-1 \choose n} C_{n-1}
\end{align*}
This implies that $n\mid S_n$ since $C_l=\bi{2l}l/(l+1)\in\Z$ for all $l\in\N$.

Now assume that $p$ is an odd prime not dividing $b(b-2c)$. For any $l\in\N$, we clearly have
$$\bi{p-1}l\bi{p+l}l=\prod_{0<j\ls l}\l(\f{p^2}{j^2}-1\r)\eq(-1)^l\ \pmod{p^2}.$$
So, by our simplified expression for $S_p$, we get
$$S_p\eq - p \sum_{l=0}^{p-2} (2lb^2-b^2+8c^2) c^{2l}(b+2c)^{p-l-1}(b-2c)^{p-l-2}(-1)^l{2l \choose  l}\pmod{p^2}.$$
By Fermat's little theorem, $(b-2c)^{p-1}\eq1\pmod p$. Note also that
$$\f{\bi{2l}l}{(-4)^l}=\bi{-1/2}l\eq\bi{(p-1)/2}l\pmod p\quad\mbox{for all}\ l=0,1,\ldots,p-1.$$
Therefore,
\begin{align*}&\f1p\sum_{k=0}^{p-1}(8ck+4c+b)\f{T_k(b,c^2)^2}{(b-2c)^{2k}}
\\\eq&\f{S_p}p\eq-\sum_{l=0}^{(p-1)/2} (2lb^2-b^2+8c^2)\f{c^{2l}(b+2c)^{p-l-1}}{(b-2c)^{l+1}}4^l\bi{(p-1)/2}l
\\=&-2b^2\cdot\f{p-1}2\sum_{l=1}^{(p-1)/2}\f{4c^2(b+2c)^{(p-1)/2}}{(b-2c)^2}\bi{(p-3)/2}{l-1}
\\&\ \ \times(b+2c)^{(p-3)/2-(l-1)}\l(\f{4c^2}{b-2c}\r)^{l-1}
\\&+\f{b^2-8c^2}{b-2c}(b+2c)^{(p-1)/2}\sum_{l=0}^{(p-1)/2}\bi{(p-1)/2}l\l(\f{4c^2}{b-2c}\r)^l(b+2c)^{(p-1)/2-l}
\\\eq&\f{4b^2c^2}{(b-2c)^2}\l(\f{b+2c}p\r)\l(\f{4c^2}{b-2c}+b+2c\r)^{(p-3)/2}
\\&+\f{b^2-8c^2}{b-2c}\l(\f{b+2c}p\r)\l(\f{4c^2}{b-2c}+b+2c\r)^{(p-1)/2}
\\=&\l(\f{4c^2}{b-2c}+\f{b^2-8c^2}{b-2c}\r)\l(\f{b+2c}p\r)\l(\f{b^2}{b-2c}\r)^{(p-1)/2}
\\\eq&(b+2c)\l(\f{b^2-4c^2}p\r)\pmod p.
\end{align*}
This concludes the proof of Theorem 1.4. \qed

\section{Proof of Theorem~\ref{s4}}

The following known result can be easily proved by Lucas' theorem (cf. \cite{HS}).

\begin{lem} \label{lem5.1} Let $n$ be a positive integer. Then $\bi{2n-1}n$ is odd if and only if $n$ is a power of $2$.
\end{lem}

\medskip
\noindent{\it Proof of Theorem} \ref{s4}. Fix an integer $n>1$. Let
\[
F(k,l) = (3k^2+k) 16^{n-1-k}  {k \choose l}^2 {2l \choose l} {2(k-l) \choose k-l}
\]
be the summand and
\[
S_n = \sum_{k=0}^{n-1} \sum_{l=0}^k F(k,l)=\sum_{k=0}^{n-1}(3k^2+k)D_k16^{n-1-k}.
\]
We find that
$$R_1(k,l) =  \frac{ 8l(k-l) }{3 k+1}\ \ \t{and}\ \
R_2(k,l)  = \frac{l^3(3 k-2 l+3)  (2 k-2 l+1)}{(3k^2+k)(k-l+1)^2}.$$
We thus derive that
\begin{align*}
S_n & = \sum_{l=0}^{n-1} R_1(n,l)F(n,l)  \\
& = 2n^3 \sum_{l=1}^{n-1} {n-1 \choose l-1} {n-1 \choose l} {2 l-1 \choose l}{2 n-2 l-1 \choose n-l}
\\&=2n^3(n-1)\sum_{l=1}^{n-1}G(n,l){2 l-1 \choose l}{2 n-2 l-1 \choose n-l},
\end{align*}
where
\begin{align*}
G(n,l)=&\f1{n-1}\bi{n-1}{l-1}\bi{n-1}l
\\=&\begin{cases}1&\t{if}\ l=1,
\\(l-1){n-2 \choose l-2} {n-1 \choose l} - (l-2) {n-2 \choose  l-1}{n-1 \choose l-1}&\t{if}\ 1<l<n.
\end{cases}
\end{align*}
Thus $a_n = S_n/(2n^3(n-1))$ is an integer.

Now we consider the parity of $a_n$. For each $l=1,\ldots,n-1$, by Lemma \ref{lem5.1},
\[
{2 l-1 \choose l}{2 n-2 l-1 \choose n-l}
\]
is odd if and only if $l=2^s$ and $n-l = 2^t$  for some $s,t\in\N$.

If $n$ is not of the form $2^s+2^t$ with $s,t\in\N$, then
\[
{2 l-1 \choose l}{2 n-2 l-1 \choose n-l}\eq0\pmod2\quad\t{for all}\ l=1,\ldots,n-1,
\]
and hence $a_n$ is even.

Suppose that $n=2^s+2^t$ with $s < t$. As
$G(n,l) = G(n, n-l)$ for all $l=1,\ldots,n-1$, we see that
\[
a_n \equiv G(n,2^s)+G(n,2^t) =2G(n,2^s)\equiv 0 \pmod{2}.
\]

Now assume that $n=2^{s+1}$ for some $s\in\N$. Then, by Lemma \ref{lem5.1} we have
\begin{align*}
a_n  \equiv G\l(n, \f n2\r) \equiv (n-1)G\l(n,\f n2\r) = {n-1 \choose n/2}^2 \equiv 1  \pmod{2}
\end{align*}
as desired.

Finally, we show the congruence (\ref{cd}) for any prime $p>3$.
Recall that
\begin{align*}S_p=&2p^3\sum_{l=1}^{p-1}\bi{p-1}{l-1}\bi{p-1}l\bi{2l-1}l\bi{2(p-l)-1}{p-l}
\\=&\f{p^3}2\sum_{l=1}^{p-1}\bi{p-1}l^2\f{l}{p-l}\bi{2l}l\bi{2(p-l)}{p-l}.
\end{align*}
By \cite[Lemma 2.1]{S11a},
$$l\bi{2l}l\bi{2(p-l)}{p-l}\eq\begin{cases}-2p\pmod{p^2}&\t{if}\ 1\ls l\ls (p-1)/2,
\\2p\pmod{p^2}&\t{if}\ (p+1)/2\ls l\ls p-1.\end{cases}$$
Clearly, $\bi{p-1}l\eq(-1)^l\pmod p$.
Thus
\begin{align*}S_p\eq&\f{p^3}2\(\sum_{l=1}^{(p-1)/2}\f{-2p}{p-l}+\sum_{l=(p+1)/2}^{p-1}\f{2p}{p-l}\)
\\\eq&p^4\(\sum_{l=1}^{(p-1)/2}\f 1l+\sum_{j=1}^{(p-1)/2}\f1j\)=2p^4H_{(p-1)/2}\pmod{p^5}.
\end{align*}
So we have
\begin{align*}\f1{p^4}&\sum_{k=0}^{p-1}\f{3k^2+k}{16^k}D_k=\f{S_p}{16^{p-1}p^4}\eq\f{S_p}{p^4}
\eq2H_{(p-1)/2}\pmod{p}.
\end{align*}
Combining this with (\ref{hh}), we finally obtain the desired (\ref{cd}).
This concludes the proof of Theorem~\ref{s4}. \qed

\section{Proof of Theorem \ref{s7}}

\medskip
\noindent{\it Proof of Theorem} \ref{s7}(i). Fix an integer $n>1$. Let
\[
F(k,l) = (12 k^4+25 k^3+21 k^2+6 k) (-1)^k {k \choose l}^2 {2l \choose k}
\]
be the summand and
\[
S_n = \sum_{k=0}^{n-1} \sum_{l=0}^k F(k,l)=\sum_{k=0}^{n-1}(12k^4+25k^3+21k^2+6k)(-1)^kf_k.
\]
We find that
$$R_1(k,l)   =  \frac{k(k-1)(k-l)^2  (3 k^2-4 k l-3 k-8 l-6)}{ (12 k^3+25 k^2+21 k+6) (-2 l+k-1) (-2 l-2+k)} $$
and
$$R_2(k,l)  = \frac{(k-1)l(k l-2 k+2 l)}{ 12 k^3+25 k^2+21 k+6 }.$$
Thus we derive that
\begin{align*}
S_n & = \sum_{l=0}^{n-1} R_1(n,l)F(n,l) \\
& = (-1)^{n-1} n^3 \sum_{l=0}^{n-1} (-3 n^2+3 n+4 n l+6+8 l) {n-1 \choose l}^2 {2l \choose n-2}.
\end{align*}
Note that
\begin{align*}&\sum_{l=0}^{n-1} (-3 n^2+3 n+4 n l+6+8 l) {n-1 \choose l}^2 {2l \choose n-2}
\\\eq&\sum_{l=0}^{n-1}(n(n-1)+6(2l+1)){n-1 \choose l}^2 {2l \choose n-2}
\\=&n(n-1)\sum_{l=0}^{n-1}\bi{n-1}l^2\bi{2l}{n-2}
\\&+6(n-1)\sum_{l=0}^{n-1}\bi{n-1}l^2\l(\bi{2l}{n-1}+\bi{2l}{n-2}\r)
\\\eq&(n-1)(n-2)\sum_{l=0}^{n-1}\bi{n-1}l^2\bi{2l}{n-2}\pmod{4(n-1)}
\end{align*}
since
$$\bi{n-1}l\bi{2l}{n-1}=\bi{n-2}{l-1}\f{n-1}l\cdot\f{2l}{n-1}\bi{2l-1}{n-2}\quad\t{for any}\ l\in\Z^+.$$
Therefore
$$(-1)^{n-1}S_n\eq n^3(n-1)(n-2)\sum_{l=0}^{n-1}\bi{n-1}l^2\bi{2l}{n-2}\pmod{4n^3(n-1)}.$$

We claim that $4\mid(n-2)\sigma_n$, where $\sigma_n=\sum_{l=0}^{n-1}\bi{n-1}l^2\bi{2l}{n-2}$.
This is trivial if $n\eq2\pmod 4$. If $n$ is odd, then
$$(n-2)\sigma_n=\sum_{l=1}^{n-1}\bi{n-1}l\bi{n-2}{l-1}\f{n-1}l\times2l\bi{2l-1}{n-3}\eq0\pmod 4.$$
When $4\mid n$, we have
\begin{align*}\sigma_n\eq&\sum_{l=0}^{n-1}\bi{n-1}l\bi{2l+1}{n-1}\bi{2l}{n-2}
\\=&\sum_{l=0}^{n-1}\bi{2l+1}{l+1}\bi{l+1}{n-1-l}=\sum_{l=1}^n\bi{2l-1}l\bi l{n-l}
\pmod 2.
\end{align*}
By Lemma \ref{lem5.1}, $\bi{2l-1}l$ with $l\in\Z^+$ is odd if and only if $l$ is a power of $2$.
If $l$ is a power of $2$ with $n/2<l<n$, then
$$\bi l{n-l}=\f{l}{n-l}\bi{l-1}{n-l-1}\eq0\pmod 2.$$
If $n$ is a power of $2$, then
$$\bi{n/2}{n-n/2}+\bi n{n-n}=2\eq0\pmod 2.$$
Therefore, when $4\mid n$ we also have $4\mid(n-2)\sigma_n$ since $\sigma_n$ is even.
So the claim always holds and hence $4n^3(n-1)\mid S_n$ as desired.

Let $p$ be an odd prime. From the above we see that
\begin{align*}S_p=&p^3\sum_{l=0}^{p-1}(-3p^2+3p+4pl+6+8l)\bi{p-1}l^2\bi{2l}{p-2}
\\\eq&2p^3\sum_{l=0}^{p-1}(4l+3)\bi{2l}{p-2}=2p^3\sum_{l=0}^{p-1}(4l+3)\bi{2l+2}p\f{p(p-1)}{(2l+1)(2l+2)}
\\\eq&2p^3\l(4\times\f{p-1}2+3\r)\bi{p+1}p\f{p(p-1)}{p(p+1)}
\\&+2p^3(4(p-1)+3)\bi{2p}p\f{p(p-1)}{(2p-1)2p}
\\\eq&2p^3(-1-1)=-4p^3\pmod{p^4}.
\end{align*}
This proves (\ref{f4}). \qed

\begin{lem} \label{lem6.1} We have
$$\sum_{k=0}^n(3k-2n)\bi nk^2\bi{2k}k=0\quad \ \mbox{for all}\ n\in\N.$$
\end{lem}
\noindent{\it Remark}. This is a known result, see, e.g., \cite[p.\,132]{PWZ}.

\begin{lem} \label{lem6.2} We have $3\mid ng_n$ for all $n\in\N$.
\end{lem}
\Proof. It is easy to see that $3\mid ng_n$ for $n=0,1$.

Let $n\in\N$ with $ng_n\eq (n+1)g_{n+1}\eq0\pmod 3$. By the Zeilberger algorithm (cf. \cite[pp.\, 101-119]{PWZ}),
$$9(n+1)^2g_n-(10n^2+30n+23)g_{n+1}+(n+2)^2g_{n+2}=0.$$
Since
$$(10n^2+30n+23)g_{n+1}\eq (n-1)(n+1)g_{n+1}\eq0\pmod 3,$$
we must have $(n+2)^2g_{n+2}\eq0\pmod 3$ and hence $3\mid(n+2)g_{n+2}$.

So far we have proved the desired result by induction. \qed

\medskip
\noindent{\it Proof of Theorem} \ref{s7}(ii). Fix $n\in\Z^+$. Let
\[
F(k,l) = (12 k^3+34 k^2+30 k+9)  {k \choose l}^2 {2l \choose l}
\]
be the summand and
\[
S_n = \sum_{k=0}^{n-1} \sum_{l=0}^k F(k,l)=\sum_{k=0}^{n-1}(12 k^3+34 k^2+30 k+9)g_k.
\]
We find the rational functions
$$R_1(k,l) =  \frac{ k (k-1)(k-l)^2 (3 k^2-4 k l-3 k-8 l-6)}{ (12 k^3+25 k^2+21 k+6) (-2 l+k-1) (-2 l-2+k)}$$
and
$$R_2(k,l)= \frac{l (k-1) (k l-2 k+2 l)}{ 12 k^3+25 k^2+21 k+6 }.$$
Thus we derive that
\begin{align*}
S_n  =& \sum_{l=0}^{n-1} R_1(n,l)F(n,l) \\
 = &n^2 \sum_{l=0}^{n-1} (-2-5 l+10 n-3 l^2+16 n l+3 n l^2) {n-1 \choose l}^2 \frac{1}{l+1} {2l \choose l}
\\&+n^4\sum_{l=0}^{n-1}(1-2l){n-1 \choose l}^2 \frac{1}{l+1} {2l \choose l}
\\\eq& n^2\sum_{l=0}^{n-1}\l(-(3l+2)(l+1)+n(l+1)+n^2(l+1)\r)\bi{n-1}l^2\f1{l+1}\bi{2l}l
\\=&n^2\sum_{l=0}^{n-1}(n^2+n-3l-2)\bi{n-1}l^2\bi{2l}l\pmod{3n^3}.
\end{align*}
By Lemma \ref{lem6.1},
$$\sum_{l=0}^{n-1}3l\bi{n-1}l^2\bi{2l}l=2(n-1)\sum_{l=0}^{n-1}\bi{n-1}l^2\bi{2l}l.$$
So we have
\begin{align*}S_n\eq&n^2\sum_{l=0}^{n-1}(n^2+n-2n)\bi{n-1}l^2\bi{2l}l=n^3(n-1)g_{n-1}\pmod{3n^3}.
\end{align*}
As $3\mid(n-1)g_{n-1}$ by Lemma \ref{lem6.2}, we obtain $S_n\eq0\pmod{3n^3}$.

Now let $p$ be an odd prime. Recall that
$$S_p=p^2\sum_{l=0}^{p-1}(-(3l+2)(l+1)+p(3l^2+16l+10)+p^2(1-2l))\bi{p-1}l^2\f1{l+1}\bi{2l}l.$$
Thus, with the help of Lemma \ref{lem6.1}, we have
\begin{align*}S_p\eq&-p^2\sum_{l=0}^{p-1}(3l+2)\bi{p-1}l^2\bi{2l}l+p^3\sum_{l=0}^{p-1}(3l^2+16l+10)\f1{l+1}\bi{2l}l
\\\eq&-2p^3\sum_{l=0}^{p-1}\bi{2l}l+p^3\sum_{l=0}^{p-1}(3l+13)\bi{2l}l-3p^3\sum_{l=0}^{p-1}C_l
\pmod{p^4}.\end{align*}
By \cite[Theorem 1.2]{PS},
$$\sum_{l=0}^{p-1}\bi{2l}l=\l(\f p3\r)\pmod p,\ \ 3\sum_{l=0}^{p-1}l\bi{2l}l\eq-2\l(\f p3\r)\pmod p,$$
and
$$\sum_{l=0}^{p-1}C_l=\sum_{l=0}^{p-1}\bi{2l}l-\sum_{l=0}^{p-1}\bi{2l}{l+1}\eq\l(\f p3\r)-\l(\f{p-1}3\r)\eq\f{3(\f p3)-1}2\pmod p.$$
Therefore, $S_p$ is congruent to
$$p^3\l(-2\l(\f p3\r)+11\l(\f p3\r)-\f 32\l(3\l(\f p3\r)-1\r)\r)=\f{3p^3}2\l(1+3\l(\f p3\r)\r)$$
modulo $p^4$. This proves (\ref{g4}). \qed

\medskip
\noindent{\it Proof of Theorem} \ref{s7}(iii). Fix $n\in\Z^+$.  Let
\[
F(k,l) = (18 k^5+45 k^4+46 k^3+24 k^2+7 k+1) (-1)^k {k+l \choose 2l}^2{2l \choose l}^2
\]
be the summand and
\[
S_n = \sum_{k=0}^{n-1} \sum_{l=0}^k F(k,l)=\sum_{k=0}^{n-1}(18 k^5+45 k^4+46 k^3+24 k^2+7 k+1) (-1)^kA_k.
\]
We find the rational functions
$$R_1(k,l)  = -\frac{k(k-l)^2  (6  k^2 l+2 k^2-6 l^2-7 l-1)}{(2 k+1) (9 k^4+18 k^3+14 k^2+5 k+1) (l+1)} $$
and
$$R_2(k,l)= \frac{l^3(3 l+4) }{ 9 k^4+18 k^3+14 k^2+5 k+1 }.$$
Thus we derive that
\begin{align*}
S_n & = \sum_{l=0}^{n-1} R_1(n,l)F(n,l) \\
& = (-1)^n n \sum_{l=0}^{n-1} (-6 n^2 l-2 n^2+6 l^2+7 l+1) {n+l \choose 2l+1}^2 \f{(2l+1)^2}{l+1} {2l \choose l}^2
\\&=(-1)^nn\sum_{l=0}^{n-1}((-6n^2+6l+1)(l+1)+4n^2)\f{n^2}{l+1}\bi{n-1}l^2\bi{n+l}l^2
\\&=(-1)^nn^3\sum_{l=0}^{n-1}(-6n^2+6l+1)\bi{n-1}l^2\bi{n+l}l^2
\\&\ \ +(-1)^n4n^4\sum_{l=0}^{n-1}\bi n{l+1}\bi{n-1}l\bi{n+l}l^2.
\end{align*}
Thus $n^4\mid S_n$ if and only if
$$\sum_{l=0}^{n-1}(6l+1)\bi{n-1}l^2\bi{n+l}l^2\eq0\pmod n.$$
Note that
\begin{align*}&\sum_{l=0}^{n-1}(6l+1)\bi{n-1}l^2\bi{n+l}l^2
\\=&3\sum_{l=0}^{n-1}(2l+1)\bi{n-1}l^2\bi{-n-1}l^2-2\sum_{l=0}^{n-1}\bi{n-1}l^2\bi{-n-1}l^2.
\end{align*}
By Guo and Zeng \cite[(1.9)]{GZ}, we have
$$\sum_{l=0}^{n-1}\bi{n-1}l^2\bi{-n-1}l^2\eq0\pmod n.$$
As proved by Sun \cite{S14b}, we also have
$$\sum_{l=0}^{n-1}(2l+1)\bi{n-1}l^2\bi{-n-1}l^2\eq0\pmod {n^2}.$$
Therefore, $S_n\eq0\pmod{n^4}$ as desired.

Let $p>3$ be a prime. Recall that
$$\f{(-1)^pS_p}{p^3}=\sum_{l=0}^{p-1}\l(-6p^2+6l+1+\f{4p^2}{l+1}\r)\bi{p-1}l^2\bi{p+l}l^2.$$
For each $l=0,\ldots,p-1$, obviously
\begin{align*}\bi{p-1}l^2\bi{p+l}l^2=&\(\prod_{0<j\ls l}\f{p^2-j^2}{j^2}\)^2
\\\eq&\(1-p^2\sum_{0<j\ls l}\f1{j^2}\)^2\eq1-2p^2\sum_{0<j\ls l}\f1{j^2}\pmod{p^4}.
\end{align*}
Thus
\begin{align*}\f{-S_p}{p^3}\eq&\sum_{l=0}^{p-1}\l(-6p^2+6l+1+\f{4p^2}{l+1}\r)
\\&-2p^2\sum_{j=1}^{p-1}\f1{j^2}\sum_{j\ls l<p}\l(-6p^2+6l+1+\f{4p^2}{l+1}\r)
\\\eq&-6p^3+6\times\f{p(p-1)}2+p+4p^2\l(H_{p-1}+\f1p\r)
\\&-2p^2\sum_{j=1}^{p-1}\f1{j^2}\(\sum_{l=j}^{p-1}(6l+1)+4p\)\pmod{p^4}.
\end{align*}
Combining this with Wolstenholme's classical congruences (cf. \cite{W})
$$H_{p-1}\eq0\pmod{p^2}\ \ \t{and}\ \ \sum_{j=1}^{p-1}\f1{j^2}\eq0\pmod p,$$
we finally obtain that
\begin{align*}
\f{-S_p}{p^3}\eq&-6p^3+3p^2-3p+p+4p-2p^2\sum_{j=1}^{p-1}\f{3(p^2-j^2)-2(p-j)}{j^2}
\\\eq&-6p^3+3p^2+2p+2p^2\times3(p-1)=2p-3p^2
\pmod{p^4},
\end{align*}
which proves (\ref{a7}).

The proof of Theorem \ref{s7} is now complete. \qed


\begin{thebibliography}{99}
\bibitem{Chen06}
W.Y.C. Chen, Q.-H. Hou, and Y.-P. Mu, {\it A telescoping method for double summations}, J. Comput. Appl. Math. {\bf 196} (2006) 553--566.

\bibitem{Fra95}
J. Franel, {\it On a question of Laisant}, L'Interm\'ediaire des Math\'ematiciens {\bf 1} (1894), 45--47.

\bibitem{Guo}
V.J.W. Guo, {\it Proof of two conjectures of Sun on congruences for Franel numbers}, Integral Transforms Spec. Funct. {\bf 24} (2013), 532--539.

\bibitem{Guo15}
V.J.W. Guo, G.-S. Mao and H. Pan, {\it Proof of a conjecture involving Sun polynomials}, J. Difference Equ. Appl. {\bf 22} (2016), 1184--1197.

\bibitem{GZa} V.J.W. Guo and J. Zeng, {\it Proof of some
conjectures of Z.-W. Sun on congruences for Ap\'ery polynomials},
J. Number Theory {\bf 132} (2012), 1731--1740.

\bibitem{GZ}
V.J.W. Guo and J. Zeng, {\it New congruences for sums involving Ap\'ery numbers or central Delannoy numbers},
Int. J. Number Theory {\bf 8} (2012), 2003--2016.


\bibitem{HS} H. Hu and Z.-W. Sun, {\it An extension of Lucas' theorem},
 ¡¡ Proc. Amer. Math. Soc. {\bf  129} (2001), 3471--3478.

\bibitem{Kou10}
C. Koutschan, {\it A fast approach to creative telescoping}, Math. Comput. Sci. {\bf 4} (2010), 259--266.

\bibitem{Lyo02}
R. Lyons, P. Paule, and A. Riese, {\it A computer proof of a series evaluation in terms of harmonic numbers}, Appl. Algebra Engrg. Comm. Comput. {\bf 13} (2002), 327--333.

\bibitem{PS}
H. Pan and Z.-W. Sun, {\it A combinatorial identity with appplications to Catalan numbers}, Discrete Math. {\bf 306} (2006), 1921--1940.

\bibitem{PWZ}
M. Petkov\v sek, H. S. Wilf and D. Zeilberger, $A=B$, A K Peters, Wellesley, 1996.

\bibitem{RS}
L.B. Richmond and J. Shallit, {\it Counting abelian squares}, Electron. J. Combin. {\bf 16} (2009), \#R72, 9 pages.

\bibitem{S} N.J.A. Sloane, {\rm Sequence A002895 in OEIS
(On-Line Encyclopedia of Integer Sequences)}, {\tt http://oeis.org/A002895}.

\bibitem{St} V. Strehl, {\it Binomial identities¨C-combinatorial and algorithmic aspects}, Discrete Math. {\bf 136} (1994), 309--346.

\bibitem{S11a} Z.-W. Sun, {\it Super congruences and Euler numbers}, Sci. China Math. {\bf 54} (2011), 2509--2535.
\bibitem{S11b} Z.-W. Sun, {\it On Delannoy numbers and Schr\"oder numbers}, J. Number Theory {\bf 131} (2011), 2387--2397.

\bibitem{S12} Z.-W. Sun, {\it On sums of Ap\'ery polynomials and
related congruences}, J. Number Theory {\bf 132} (2012), 2673--2699.

\bibitem{S13a} Z.-W. Sun, {\it Congruences for Franel numbers},
 Adv. in Appl. Math. {\bf 51} (2013), 524--535.

\bibitem{S13b} Z.-W. Sun, {\it Connections between $p=x^2+3y^2$ and Franel numbers},
 ¡¡¡¡J. Number Theory {\bf 133} (2013), 2914--2928.

\bibitem{S13c} Z.-W. Sun,  Conjectures and results on $x^2$ mod $p^2$ with $4p=x^2+dy^2$,
in: Number Theory and Related Area (eds., Y. Ouyang, C. Xing, F. Xu and P. Zhang),
¡¡¡¡Adv. Lect. Math. 27, Higher Education Press and International Press, Beijing-Boston, 2013, pp. 149-197.


\bibitem{S14} Z.-W. Sun, {\it Congruences involving generalized central trinomial coefficients}, Sci. China Math. {\bf 57} (2014), 1375--1400.

\bibitem{S14b} Z.-W. Sun, {\it Two new kinds of numbers and related divisibility results}, {\tt arXiv:1408.5381} [math.NT], preprint, 2014.

\bibitem{Sun16}
Z.-W. Sun, {\it Congruences involving $g_n(x) = \sum_{k=0}^n {n \choose k}^2 {2k \choose k} x^k$}, Ramanujan J. {\bf 40} (2016), 511--533.

\bibitem{W}
J. Wolstenholme, {\it On certain properties of prime numbers}, Quart. J. Appl. Math.
{\bf 5} (1862), 35--39.
\end{thebibliography}
\end{document}